\documentstyle[amssymb,amsmath]{article}
\newtheorem{th}{Theorem}[section]
\newtheorem{lem}[th]{Lemma}

\newtheorem{cor}[th]{Corollary}
\newtheorem{defn}[th]{Definition}
\newenvironment{defn-new}{\begin{defn} \em}{\end{defn}}
\newtheorem{rem}[th]{Remark}
\newenvironment{rem-new}{\begin{rem} \em}{\end{rem}}
\newtheorem{ex}[th]{Example}
\newenvironment{ex-new}{\begin{ex} \em}{\end{ex}}

\newenvironment{notation-new}{\begin{rem} \em}{\end{rem}}

\newenvironment{agr-new}{\begin{rem} \em}{\end{rem}}

\makeatletter \@addtoreset{equation}{section} \makeatother

\makeatletter \@addtoreset{figure}{section} \makeatother

\begin{document}

\begin{center}
{\Large {\bf On $3$-dimensional $\left(\varepsilon \right)$-para Sasakian
manifold}}

\bigskip \bigskip {\large {\bf Punam Gupta}}
\end{center}

\noindent {\bf Abstract.} The purpose of the present paper is to study the
globally and locally $\varphi $-${\cal T}$-symmetric $\left( \varepsilon
\right) $-para Sasakian manifold in dimension $3$. The globally $\varphi $-$%
{\cal T}$-symmetric $3$-dimensional $\left( \varepsilon \right) $-para
Sasakian manifold is either Einstein manifold or has a constant scalar
curvature. The necessary and sufficient condition for Einstein manifold to
be globally $\varphi $-${\cal T}$ -symmetric is given. A $3$-dimensional $%
\left( \varepsilon \right) $ -para Sasakian manifold is locally $\varphi $-$%
{\cal T}$-symmetric if and only if the scalar curvature $r$ is constant. A $%
3 $-dimensional $\left( \varepsilon \right) $-para Sasakian manifold with $%
\eta $-parallel Ricci tensor is locally $\varphi $-${\cal T}$-symmetric. In
the last, an example of $3$-dimensional locally $\varphi $-${\cal T}$%
-symmetric $\left( \varepsilon \right) $-para Sasakian manifold is given.
\medskip

 \medskip

\noindent {\bf Keywords.} ${\cal T}$-curvature tensor; $\left( \varepsilon
\right) $-para Sasakian manifold; globally and locally $\varphi $-${\cal T}$%
-symmetric manifold; $\eta $-parallel Ricci tensor. \medskip

\section{Introduction}

Let $M$ be an $m$-dimensional semi-Riemannian manifold and $\nabla $ the
Levi-Civita connection on $M$. A semi-Riemannian manifold $M$ is said to
recurrent \cite{Ruse-49} if the Riemann curvature tensor $R$ satisfies the
relation 
\[
\left( \nabla _{U}R\right) \left( X,Y,Z,V\right) =\alpha (U)R\left(
X,Y,Z,V\right) ,\quad X,Y,Z,V,U\in TM, 
\]%
where $\alpha $ is $1$-form. If $\alpha =0$, then $M$ is called symmetric in
the sense of Cartan \cite{Cartan-26}. \medskip

In 1977, Takahashi \cite{Takahashi-77} introduced the notion of locally $%
\varphi $-symmetry on a Sasakian manifold, which is weaker than the local
symmetry. A Sasakian manifold is said to have locally $\varphi $-symmetry if
it satisfies 
\[
\varphi ^{2}\left( \left( \nabla _{U}R\right) (X,Y)Z\right) =0, 
\]%
where $X$, $Y$, $Z$, $U$ are horizontal vector fields. If $X$, $Y$, $Z$, $U$
are arbitrary vector fields, then it is known as globally $\varphi $%
-symmetric Sasakian manifold. A $\varphi $-symmetric space condition is weak
condition for a Sasakian manifold in comparision to the symmetric space
condition. Local symmetry is a very strong condition for the class of $K$%
-contact or Sasakian manifolds. Indeed, such spaces must have constant
curvature equal to $1$ (\cite{Okumura-62}, \cite{Tanno-68}). On the other
hand, local symmetry is also a very strong condition for the class of $%
\left( \varepsilon \right) $-para Sasakian manifold. Such spaces must have
constant curvature equal to $-\,\varepsilon $ \cite{TKYK-09}. In 2010,
Tripathi et al. \cite{TKYK-09} proved that the condition of semi-symmetry ($%
R\cdot R=0$), symmetry and have a constant curvature $-\,\varepsilon $ is
equivalent for $\left( \varepsilon \right) $-para Sasakian manifold. \medskip

Three-dimensional locally $\varphi $-symmetric Sasakian manifold is studied
by Watanabe \cite{Watanabe}. Many authors like De \cite{De-08}, De et al. 
\cite{De-09}, De and Pathak \cite{De-Pathak-04}, Shaikh and De \cite%
{Shaikh-De-00} have extended this notion to $3$-dimensional Kenmotsu,
trans-Sasakian and LP-Sasakian manifolds. Yildiz et al. \cite{Yildiz}
studied the case for $3$-dimensional $\alpha $-Sasakian manifolds and gave
the example for locally $\varphi $-symmetric $3$-dimensional $\alpha $%
-Sasakian manifolds. De and De \cite{De-De} studied the $\varphi $%
-concircularly symmetric Kenmotsu manifold and gave the example of such
manifold in dimension $3$. De et al. \cite{De-Ozgur-09} studied the $3$%
-dimensional globally and locallly $\varphi $-quasiconformally symmetric
Sasakian manifolds and also gave the example. \medskip

In the present work, globally and locally $\varphi $-${\cal T}$-symmetric $%
\left( \varepsilon \right) $-para Sasakian manifold in dimension $3$ is
studied. The paper is organized as follows: Section $2$ and $3$ is devoted
to the study of ${\cal T}$-curvature tensor and $\left( \varepsilon \right) $%
-para Sasakian manifold, respectively. Some results for $3$-dimensional $%
\left( \varepsilon \right) $-para Sasakian manifold are given. The necessary
and sufficient condition for $\left( \varepsilon \right) $-para Sasakian
manifold of constant curvature is given. In section $4$, the definition of
globally and locally $\varphi $-${\cal T}$-symmetric manifold are given.
Globally $\varphi $-${\cal T}$-symmetric $3$-dimensional $\left( \varepsilon
\right) $-para Sasakian manifold is either Einstein or has a constant scalar
curvature under some condition. The necessary and sufficient condition for
locally $\varphi $-${\cal T}$-symmetric $3$-dimensional $\left( \varepsilon
\right) $-para Sasakian manifold to be locally $\varphi $-symmetric is
given. In section $5$, the definition of $\eta $-parallel $\left(
\varepsilon \right) $-para Sasakian manifold is given. A $3$-dimensional $%
\left( \varepsilon \right) $-para Sasakian manifold with $\eta $-parallel
Ricci tensor is locally $\varphi $-${\cal T}$-symmetric. In the last
section, the example of a locally $\varphi $-${\cal T}$-symmetric $3$%
-dimensional $\left( \varepsilon \right) $-para Sasakian manifold is given.

\section{${\cal T}$-curvature tensor\label{sect-TCT}}

The definition of ${\cal T}$-curvature tensor \cite{Tripathi-Gupta}\ is
given by

\begin{defn-new}
\label{defn-GCT} In an $m$-dimensional semi-Riemannian manifold $\left(
M,g\right) $, the {\em ${\cal T}$-curvature tensor} of type $(1,3)$ defined
by 
\begin{eqnarray}
{\cal T}\left( X,Y\right) Z &=&a_{0}\,R\left( X,Y\right) Z  \nonumber \\
&&+\ a_{1}\,S\left( Y,Z\right) X+a_{2}\,S\left( X,Z\right) Y+a_{3}\,S(X,Y)Z 
\nonumber \\
&&+\ a_{4}\,g\left( Y,Z\right) QX+a_{5}\,g\left( X,Z\right)
QY+a_{6}\,g(X,Y)QZ  \nonumber \\
&&+\ a_{7}\,r\left( g\left( Y,Z\right) X-g\left( X,Z\right) Y\right) ,
\label{eq-GCT-1}
\end{eqnarray}%
for all $X,Y,Z\in TM$, where $a_{0},\ldots ,a_{7}$ are some constants; and $%
R $, $S$, $Q$ and $r$ are the curvature tensor, the Ricci tensor, the Ricci
operator of type $(1,1)$ and the scalar curvature respectively.
\end{defn-new}

In particular, the ${\cal T}$-curvature tensor is reduced to

\begin{enumerate}
\item the {\em Riemann curvature tensor} $R$ if 
\[
a_{0}=1,\quad a_{1}=a_{2}=a_{3}=a_{4}=a_{5}=a_{6}=a_{7}=0, 
\]

\item the {\em quasiconformal curvature tensor} ${\cal C}_{\ast }$ \cite%
{Yano-Sawaki-68} if 
\[
a_{1}=-a_{2}=a_{4}=-a_{5},\quad a_{3}=a_{6}=0,\quad a_{7}=-\,\frac{1}{m}%
\left( \frac{a_{0}}{m-1}+2a_{1}\right) , 
\]

\item the {\em conformal curvature tensor} ${\cal C}$ \cite[p.~90]%
{Eisenhart-49} if 
\[
a_{0}=1,\quad a_{1}=-a_{2}=a_{4}=-a_{5}=-\,\frac{1}{m-2},\quad
a_{3}=a_{6}=0,\quad a_{7}=\frac{1}{(m-1)(m-2)}, 
\]

\item the {\em conharmonic curvature tensor} ${\cal L}$ \cite{Ishii-57} if 
\[
a_{0}=1,\quad a_{1}=-a_{2}=a_{4}=-a_{5}=-\,\frac{1}{m-2},\,\quad
a_{3}=a_{6}=0,\quad a_{7}=0, 
\]

\item the {\em concircular curvature tensor} ${\cal V}$ (\cite{Yano-40}, 
\cite[p. 87]{Yano-Bochner-53}) if 
\[
a_{0}=1,\quad a_{1}=a_{2}=a_{3}=a_{4}=a_{5}=a_{6}=0,\quad a_{7}=-\,\frac{1}{%
m(m-1)}, 
\]

\item the {\em pseudo-projective curvature tensor }${\cal P}_{\ast }$ \cite%
{Prasad-2002} if 
\[
a_{1}=-a_{2},\quad a_{3}=a_{4}=a_{5}=a_{6}=0,\quad a_{7}=-\,\frac{1}{m}%
\left( \frac{a_{0}}{m-1}+a_{1}\right) , 
\]

\item the {\em projective curvature tensor} ${\cal P}$ \cite[p. 84]%
{Yano-Bochner-53} if 
\[
a_{0}=1,\quad a_{1}=-a_{2}=-\,\frac{1}{(m-1)}\text{,\quad }%
a_{3}=a_{4}=a_{5}=a_{6}=a_{7}=0, 
\]

\item the $M${\em -projective curvature tensor }\cite{Pokhariyal-Mishra-71}
if 
\[
a_{0}=1,\quad a_{1}=-a_{2}=a_{4}=-a_{5}=-\frac{1}{2(m-1)},\quad
a_{3}=a_{6}=a_{7}=0, 
\]

\item the $W_{0}$-{\em curvature tensor} \cite[eq (1.4)]%
{Pokhariyal-Mishra-71} if 
\[
a_{0}=1,\quad a_{1}=-a_{5}=-\,\frac{1}{(m-1)},\quad
a_{2}=a_{3}=a_{4}=a_{6}=a_{7}=0, 
\]

\item the $W_{0}^{\ast }$-{\em curvature tensor} \cite[eq (1.4)]%
{Pokhariyal-Mishra-71} if 
\[
a_{0}=1,\quad a_{1}=-a_{5}=\frac{1}{(m-1)},\quad
a_{2}=a_{3}=a_{4}=a_{6}=a_{7}=0, 
\]

\item the $W_{1}$-{\em curvature tensor} \cite{Pokhariyal-Mishra-71} if 
\[
a_{0}=1,\quad a_{1}=-a_{2}=\frac{1}{(m-1)},\quad
a_{3}=a_{4}=a_{5}=a_{6}=a_{7}=0, 
\]

\item the $W_{1}^{\ast }$-{\em curvature tensor} \cite{Pokhariyal-Mishra-71}
if 
\[
a_{0}=1,\quad a_{1}=-a_{2}=-\,\frac{1}{(m-1)},\quad
a_{3}=a_{4}=a_{5}=a_{6}=a_{7}=0, 
\]

\item the $W_{2}$-{\em curvature tensor} \cite{Pokhariyal-Mishra-70} if 
\[
a_{0}=1,\quad a_{4}=-a_{5}=-\,\frac{1}{(m-1)},\quad
a_{1}=a_{2}=a_{3}=a_{6}=a_{7}=0, 
\]

\item the $W_{3}$-{\em curvature tensor} \cite{Pokhariyal-Mishra-71} if 
\[
a_{0}=1,\quad a_{2}=-a_{4}=-\,\frac{1}{(m-1)},\quad
a_{1}=a_{3}=a_{5}=a_{6}=a_{7}=0, 
\]

\item the $W_{4}$-{\em curvature tensor} \cite{Pokhariyal-Mishra-71} if 
\[
a_{0}=1,\quad a_{5}=-a_{6}=\frac{1}{(m-1)},\quad
a_{1}=a_{2}=a_{3}=a_{4}=a_{7}=0, 
\]

\item the $W_{5}$-{\em curvature tensor} \cite{Pokhariyal-82} if 
\[
a_{0}=1,\quad a_{2}=-a_{5}=-\,\frac{1}{(m-1)},\quad
a_{1}=a_{3}=a_{4}=a_{6}=a_{7}=0, 
\]

\item the $W_{6}$-{\em curvature tensor} \cite{Pokhariyal-82} if 
\[
a_{0}=1,\quad a_{1}=-a_{6}=-\,\frac{1}{(m-1)},\quad
a_{2}=a_{3}=a_{4}=a_{5}=a_{7}=0, 
\]

\item the $W_{7}$-{\em curvature tensor} \cite{Pokhariyal-82} if 
\[
a_{0}=1,\quad a_{1}=-a_{4}=-\,\frac{1}{(m-1)},\quad
a_{2}=a_{3}=a_{5}=a_{6}=a_{7}=0, 
\]

\item the $W_{8}$-{\em curvature tensor} \cite{Pokhariyal-82} if 
\[
a_{0}=1,\quad a_{1}=-a_{3}=-\,\frac{1}{(m-1)},\quad
a_{2}=a_{4}=a_{5}=a_{6}=a_{7}=0, 
\]

\item the $W_{9}$-{\em curvature tensor} \cite{Pokhariyal-82} if 
\[
a_{0}=1,\quad a_{3}=-\,a_{4}=\frac{1}{(m-1)},\quad
a_{1}=a_{2}=a_{5}=a_{6}=a_{7}=0. 
\]
\end{enumerate}

\section{$\left( \protect\varepsilon \right) $-para Sasakian manifold}

A manifold $M$ is said to admit an almost paracontact structure if it admit
a tensor field $\varphi $ of type $\left( 1,1\right) $, a vector field $\xi $
and a $1$-form $\eta $ satisfying 
\begin{equation}
\varphi ^{2}=I-\eta \otimes \xi ,\quad \eta \left( \xi \right) =1,\quad
\varphi \xi =0,\quad \eta \circ \varphi =0.  \label{eq-1-1}
\end{equation}%
Let $g$ be a semi-Riemannian metric with ${\rm index}(g)=\nu $ such that 
\begin{equation}
g\left( \varphi X,\varphi Y\right) =g\left( X,Y\right) -\varepsilon \eta
(X)\eta \left( Y\right) ,\qquad X,Y\in TM,  \label{eq-metric-1}
\end{equation}%
where $\varepsilon =\pm 1$. Then $M$ is called an $\left( \varepsilon
\right) $-{\em almost paracontact metric manifold} equipped with an $\left(
\varepsilon \right) ${\em -almost paracontact metric structure} $(\varphi
,\xi ,\eta ,g,\varepsilon )$. In particular, if ${\rm index}(g)=1$, then an $%
(\varepsilon )$-almost paracontact metric manifold is said to be a {\em %
Lorentzian almost paracontact manifold}. In particular, if the metric $g$ is
positive definite, then an $(\varepsilon )$-almost paracontact metric
manifold is the usual {\em almost paracontact metric manifold} \cite{Sato-76}%
.

The equation (\ref{eq-metric-1}) is equivalent to 
\begin{equation}
g\left( X,\varphi Y\right) =g\left( \varphi X,Y\right)  \label{eq-metric-2}
\end{equation}%
along with 
\begin{equation}
g\left( X,\xi \right) =\varepsilon \eta (X).  \label{eq-metric-3}
\end{equation}%
From (\ref{eq-1-1}) and (\ref{eq-metric-3}) it follows that 
\begin{equation}
g\left( \xi ,\xi \right) =\varepsilon .  \label{eq-g(xi,xi)}
\end{equation}

\begin{defn-new}
\label{def-eps-pS} An $\left( \varepsilon \right) $-almost paracontact
metric structure is called an $\left( \varepsilon \right) ${\em -para
Sasakian structure} if 
\begin{equation}
\left( \nabla _{X}\varphi \right) Y=-\,g(\varphi X,\varphi Y)\xi
-\varepsilon \eta \left( Y\right) \varphi ^{2}X,\qquad X,Y\in TM,
\label{eq-eps-PS-def-1}
\end{equation}%
where $\nabla $ is the Levi-Civita connection with respect to $g$. A
manifold endowed with an $\left( \varepsilon \right) $-para Sasakian
structure is called an $\left( \varepsilon \right) ${\em -para Sasakian
manifold} \cite{TKYK-09}.
\end{defn-new}

For $\varepsilon =1$ and $g$ Riemannian, $M$ is the usual para Sasakian
manifold \cite{Sato-77,Sasaki-80}. For $\varepsilon =-1$, $g$ Lorentzian and 
$\xi $ replaced by $-\xi $, $M$ becomes a Lorentzian para Sasakian manifold 
\cite{Mat-89}.

For $\left( \varepsilon \right) $-para Sasakian manifold, it is easy to
prove that 
\begin{equation}
R(X,Y)\xi =\eta (X)Y-\eta (Y)X,  \label{eq-curvature}
\end{equation}%
\begin{equation}
R(\xi ,X)Y=\eta (Y)X-\varepsilon g(X,Y)\xi ,  \label{eq-curvature-2}
\end{equation}%
\begin{equation}
R(\xi ,X)\xi =X-\eta (X)\xi ,  \label{eq-curvature-3}
\end{equation}%
\begin{equation}
R\left( X,Y,Z,\xi \right) =\eta \left( Y\right) g\left( X,Z\right) -\eta
\left( X\right) g\left( Y,Z\right) ,  \label{eq-eps-PS-R(X,Y,Z,xi)}
\end{equation}%
\begin{equation}
\eta \left( R\left( X,Y\right) Z\right) =\varepsilon (\eta \left( Y\right)
g\left( X,Z\right) -\eta \left( X\right) g\left( Y,Z\right) ),
\label{eq-eps-PS-eta(R(X,Y),Z)}
\end{equation}%
\begin{equation}
S(X,\xi )=-(m-1)\eta (X),  \label{eq-ricci}
\end{equation}%
\begin{equation}
Q\xi =-\,\varepsilon (m-1)\xi ,  \label{eq-Q}
\end{equation}%
\begin{equation}
S(\xi ,\xi )=-(m-1),  \label{eq-S-xi-xi}
\end{equation}%
\begin{equation}
S(\varphi X,\varphi Y)=S(Y,Z)+(m-1)\eta \left( X\right) \eta \left( Y\right)
,  \label{eq-phi-phi}
\end{equation}%
\begin{equation}
\nabla _{X}\xi =\varepsilon \varphi X.  \label{eq-del}
\end{equation}%
For detail study of $\left( \varepsilon \right) $-para Sasakian manifold,
see \cite{TKYK-09}.

It is well known that in a $3$-dimensional semi-Riemannian manifold the
conformal curvature tensor ${\cal C}$ vanishes, therefore 
\begin{eqnarray}
R(X,Y)Z &=&g(Y,Z)QX-g(X,Z)QY+S(Y,Z)X-S(X,Z)Y  \nonumber \\
&&-\frac{r}{2}(g(Y,Z)X-g(X,Z)Y).  \label{eq-R3}
\end{eqnarray}

\begin{th}
Let $M$ be a $3$-dimensional $\left( \varepsilon \right) $-para Sasakian
manifold. Then 
\begin{equation}
QX=\left( \frac{r}{2}+\varepsilon \right) X-\left( \frac{r}{2}+3\varepsilon
\right) \eta (X)\xi .  \label{eq-R11}
\end{equation}
\end{th}

\noindent {\bf Proof.} Take $Z=\xi $ in (\ref{eq-R3}) and using (\ref%
{eq-metric-3}), (\ref{eq-curvature}), (\ref{eq-ricci}), we get 
\begin{equation}
\left( \frac{\varepsilon r}{2}+1\right) (\eta (Y)X-\eta (X)Y)=\varepsilon
\left( \eta (Y)QX-\eta (X)QY\right) .  \label{eq-R31}
\end{equation}%
Putting $Y=\xi $ in (\ref{eq-R31}) and using (\ref{eq-Q}), we get (\ref%
{eq-R11}).

\begin{cor}
Let $M$ be a $3$-dimensional $\left( \varepsilon \right) $-para Sasakian
manifold. Then 
\begin{equation}
S(X,Y)=\left( \frac{r}{2}+\varepsilon \right) g(X,Y)-\left( \frac{%
\varepsilon r}{2}+3\right) \eta (X)\eta (Y)  \label{eq-S-3}
\end{equation}%
and 
\begin{eqnarray}
R(X,Y)Z &=&\left( \frac{r}{2}+2\varepsilon \right) (g(Y,Z)X-g(X,Z)Y) 
\nonumber \\
&&+\left( \frac{\varepsilon r}{2}+3\right) (\eta (X)\eta (Z)Y-\eta (Y)\eta
(Z)X)  \nonumber \\
&&+\left( \frac{r}{2}+3\varepsilon \right) (g(X,Z)\eta (Y)\xi -g(Y,Z)\eta
(X)\xi ).  \label{eq-S-4}
\end{eqnarray}
\end{cor}

\begin{lem}
A $3$-dimensional $\left( \varepsilon \right) $-para Sasakian manifold is a
manifold of constant curvature if and only if $r=-6\varepsilon $.
\end{lem}

\begin{cor}
Let $M$ be a $3$-dimensional $\left( \varepsilon \right) $-para Sasakian
manifold. Then 
\begin{eqnarray}
{\cal T}(X,Y)Z &=&\left( \left( \frac{r}{2}+\varepsilon \right) \left(
a_{0}+a_{1}+a_{4}\right) +a_{7}r+\varepsilon a_{0}\right) g(Y,Z)X  \nonumber
\\
&&-\left( \left( \frac{r}{2}+\varepsilon \right) \left(
a_{0}-a_{2}-a_{5}\right) +a_{7}r+\varepsilon a_{0}\right) g(X,Z)Y  \nonumber
\\
&&+\left( \frac{r}{2}+\varepsilon \right) \left( a_{3}+a_{6}\right)
g(X,Y)Z-\left( \frac{\varepsilon r}{2}+3\right) a_{3}\eta (X)\eta (Y)Z 
\nonumber \\
&&-\left( \frac{\varepsilon r}{2}+3\right) \left( a_{0}+a_{1}\right) \eta
(Y)\eta (Z)X+\left( \frac{\varepsilon r}{2}+3\right) \left(
a_{0}-a_{2}\right) \eta (X)\eta (Z)Y  \nonumber \\
&&+\left( \frac{r}{2}+3\varepsilon \right) \left( a_{0}-a_{5}\right)
g(X,Z)\eta (Y)\xi -\left( \frac{r}{2}+3\varepsilon \right) a_{6}g(X,Y)\eta
(Z)\xi  \nonumber \\
&&-\left( \frac{r}{2}+3\varepsilon \right) \left( a_{0}+a_{4}\right)
g(Y,Z)\eta (X)\xi .  \label{eq-S-5}
\end{eqnarray}
\end{cor}

\section{$3$-dimensional $\protect\varphi $-${\cal T}$-symmetric $\left( 
\protect\varepsilon \right) $-para Sasakian manifold}

We begin with the following definition.

\begin{defn-new}
An $\left( \varepsilon \right) $-para Sasakian manifold is said to be
locally $\varphi $-${\cal T}$-symmetric manifold if 
\begin{equation}
\varphi ^{2}((\nabla _{W}{\cal T})(X,Y)Z)=0,  \label{eq-phi-T-sym}
\end{equation}%
for arbitrary vector fields $X,Y,Z,W$ orthogonal to $\xi $. If $X$, $Y$, $Z$%
, $W$ are arbitrary vector fields, then it is known as globally $\varphi $-$%
{\cal T}$-symmetric manifold.
\end{defn-new}

This notion of locally $\varphi $-symmetric was introduced by Takahashi for
Sasakian manifolds \cite{Takahashi-77}.

\begin{th}
\label{thm-T} Let $M$ be a $m$-dimensional globally $\varphi $-${\cal T}$%
-symmetric $\left( \varepsilon \right) $-para Sasakian manifold. Then

\begin{description}
\item[(i)] $M$ is Einstein manifold if $a_{0}+(m-1)a_{1}+a_{2}+a_{6}\not=0$.

\item[(ii)] $M$ has constant scalar curvature if $%
a_{0}+(m-1)a_{1}+a_{2}+a_{6}=0$ and $a_{4}+(m-1)a_{7}\not=0$.
\end{description}
\end{th}

\noindent {\bf Proof.} Let $M$ be a $m$-dimensional $\varphi $-${\cal T}$%
-symmetric $\left( \varepsilon \right) $-para Sasakian manifold. Then by
using (\ref{eq-1-1}) and (\ref{eq-phi-T-sym}), we have 
\[
(\nabla _{W}{\cal T})(X,Y)Z-\eta \left( (\nabla _{W}{\cal T})(X,Y)Z\right)
\xi =0, 
\]%
from which it follows that 
\begin{equation}
g((\nabla _{W}{\cal T})(X,Y)Z,U)-\eta \left( (\nabla _{W}{\cal T}%
)(X,Y)Z\right) g(\xi ,U)=0.  \label{eq-1-2}
\end{equation}%
Using (\ref{eq-GCT-1}) in (\ref{eq-1-2}), we obtain 
\begin{eqnarray}
0 &=&a_{0}\,(\nabla _{W}R)(X,Y,Z,U)+\ a_{1}\,(\nabla _{W}S)\left( Y,Z\right)
g(X,U)+a_{2}\,(\nabla _{W}S)\left( X,Z\right) g(Y,U)  \nonumber \\
&&+a_{3}\,(\nabla _{W}S)(X,Y)g(Z,U)+a_{4}\,(\nabla _{W}S)(X,U)g\left(
Y,Z\right) +a_{5}\,(\nabla _{W}S)(Y,U)g\left( X,Z\right)  \nonumber \\
&&+a_{6}\,(\nabla _{W}S)(Z,U)g(X,Y)+a_{7}\,(\nabla _{W}r)\left( g\left(
Y,Z\right) g(X,U)-g\left( X,Z\right) g(Y,U)\right)  \nonumber \\
&&+\eta (U)\left( a_{0}\,(\nabla _{W}R)(X,Y,Z,\xi )+\ a_{1}\,(\nabla
_{W}S)\left( Y,Z\right) g(X,\xi )+a_{2}\,(\nabla _{W}S)\left( X,Z\right)
g(Y,\xi )\right.  \nonumber \\
&&+a_{3}\,(\nabla _{W}S)(X,Y)g(Z,\xi )+\ a_{4}\,g\left( Y,Z\right) (\nabla
_{W}S)(X,\xi )+a_{5}\,g\left( X,Z\right) (\nabla _{W}S)(Y,\xi )  \nonumber \\
&&\left. +a_{6}\,g(X,Y)(\nabla _{W}S)(Z,\xi )+\ a_{7}\,(\nabla _{W}r)\left(
g\left( Y,Z\right) g(X,\xi )-g\left( X,Z\right) g(Y,\xi )\right) \right) .
\label{eq-1-3}
\end{eqnarray}%
Let $\{e_{i}\}$, $i=1,\ldots ,m$ be an orthonormal basis of tangent space at
any point of the manifold. Taking $X=U=e_{i}$ in (\ref{eq-1-3}), we get 
\begin{eqnarray}
0 &=&\left( a_{0}+(m-1)a_{1}+a_{2}+a_{3}+a_{5}+a_{6}\right) \,(\nabla
_{W}S)(Y,Z)-a_{0}\varepsilon \,\sum\limits_{i=1}^{m}(\nabla _{W}R)\left(
e_{i},Y,Z,\xi \right) g(e_{i},\xi )  \nonumber \\
&&+\left( a_{4}+(m-1)a_{7}\right) \,(\nabla _{W}r)g\left( Y,Z\right)
+a_{7}\,(\nabla _{W}r)\left( g\left( Y,Z\right) -\varepsilon \eta (Y)\eta
(Z)\right)  \nonumber \\
&&-\left( a_{2}+a_{6}\right) \,(\nabla _{W}S)(Z,\xi )\eta (Y)-\left(
a_{3}+a_{5}\right) \,(\nabla _{W}S)(Y,\xi )\eta (Z).  \label{eq-1-4}
\end{eqnarray}%
Putting $Z=\xi $ in (\ref{eq-1-4}), we have 
\begin{eqnarray}
0 &=&\left( a_{0}+(m-1)a_{1}+a_{2}+a_{6}\right) \,(\nabla _{W}S)(Y,\xi ) 
\nonumber \\
&&-a_{0}\varepsilon \,\sum\limits_{i=1}^{m}(\nabla _{W}R)\left( e_{i},Y,\xi
,\xi \right) g(e_{i},\xi )  \nonumber \\
&&+\left( a_{4}+(m-1)a_{7}\right) \,(\nabla _{W}r)g\left( Y,\xi \right) 
\nonumber \\
&&-\left( a_{2}+a_{6}\right) \,(\nabla _{W}S)(\xi ,\xi )\eta (Y).
\label{eq-1-5}
\end{eqnarray}%
Since, we have 
\begin{eqnarray}
(\nabla _{W}R)\left( e_{i},Y,\xi ,\xi \right) &=&g((\nabla _{W}R)\left(
e_{i},Y)\xi ,\xi \right)  \nonumber \\
&=&g(\nabla _{W}R\left( e_{i},Y)\xi ,\xi \right) -g(R\left( \nabla
_{W}e_{i},Y)\xi ,\xi \right)  \nonumber \\
&&-g(R\left( e_{i},\nabla _{W}Y)\xi ,\xi \right) -g(R\left( e_{i},Y)\nabla
_{W}\xi ,\xi \right)  \label{eq-1-6}
\end{eqnarray}%
at any point $p\in M$. We know that $\left\{ e_{i}\right\} $ is an
orthonormal basis, therefore $\nabla _{W}e_{i}=0$ at $p$. Using (\ref%
{eq-metric-3}) and (\ref{eq-curvature}) in (\ref{eq-1-6}), we have 
\begin{equation}
(\nabla _{W}R)\left( e_{i},Y,\xi ,\xi \right) =g(\nabla _{W}R\left(
e_{i},Y)\xi ,\xi \right) -g(R\left( e_{i},Y)\nabla _{W}\xi ,\xi \right) .
\label{eq-1-7}
\end{equation}%
By using the property of curvature tensor 
\[
g(R(e_{i},Y)\xi ,\xi )=-g(R(\xi ,\xi )Y,e_{i})=0, 
\]%
we have 
\begin{equation}
g(\nabla _{W}R(e_{i},Y)\xi ,\xi )+g(R(e_{i},Y)\xi ,\nabla _{W}\xi )=0.
\label{eq-1-8}
\end{equation}%
By (\ref{eq-1-7}) and (\ref{eq-1-8}), we get 
\begin{equation}
(\nabla _{W}R)\left( e_{i},Y,\xi ,\xi \right) =0.  \label{eq-1-9}
\end{equation}%
We know that 
\begin{equation}
(\nabla _{W}S)(Y,\xi )=\nabla _{W}S(Y,\xi )-S(\nabla _{W}Y,\xi )-S(Y,\nabla
_{W}\xi ).  \label{eq-1-10}
\end{equation}%
Using (\ref{eq-ricci}), (\ref{eq-del}) in (\ref{eq-1-10}), we get 
\begin{eqnarray}
(\nabla _{W}S)(Y,\xi ) &=&\nabla _{W}(-(m-1)\eta (Y))+(m-1)\eta \left(
\nabla _{W}Y\right) -S(Y,\varepsilon \varphi W)  \nonumber \\
&=&-(m-1)\varepsilon g(Y,\varepsilon \varphi W)-\varepsilon S(Y,\varphi W) 
\nonumber \\
&=&-(m-1)g(Y,\varphi W)-\varepsilon S(Y,\varphi W).  \label{eq-1-11}
\end{eqnarray}%
By (\ref{eq-1-11}), we have 
\begin{equation}
(\nabla _{W}S)(\xi ,\xi )=0.  \label{eq-1-12}
\end{equation}%
Using (\ref{eq-1-9}), (\ref{eq-1-11}), (\ref{eq-1-12}) in (\ref{eq-1-5}), we
have 
\begin{eqnarray}
0 &=&\left( a_{0}+(m-1)a_{1}+a_{2}+a_{6}\right) \left( -(m-1)g(Y,\varphi
W)-\varepsilon S(Y,\varphi W)\right)  \nonumber \\
&&+\varepsilon \left( a_{4}+(m-1)a_{7}\right) \,(\nabla _{W}r)\eta (Y).
\label{eq-1-13}
\end{eqnarray}%
Replacing $Y$ by $\varphi Y$ in (\ref{eq-1-13}) and using (\ref{eq-metric-1}%
), (\ref{eq-phi-phi}), we get 
\[
S(Y,W)=-\varepsilon (m-1)g(Y,W),\quad a_{0}+(m-1)a_{1}+a_{2}+a_{6}\not=0. 
\]%
If $a_{0}+(m-1)a_{1}+a_{2}+a_{6}=0$ and $a_{4}+(m-1)a_{7}\not=0$, then by (%
\ref{eq-1-5}), we have $\nabla _{W}r=0$, that is, $r=$ constant.

\begin{rem-new}
The first condition of Theorem~\ref{thm-T} is satisfied if ${\cal T}\in
\left\{ R,{\cal C}_{\ast },{\cal V},{\cal P}_{\ast },{\cal P},{\cal M},{\cal %
W}_{0}^{\ast },{\cal W}_{1},{\cal W}_{1}^{\ast },{\cal W}_{2},\ldots
,\right. \linebreak \left.{\cal W}_{6},{\cal W}_{9}\right\} $ and second
condition is satisfied if ${\cal T}\in \left\{ {\cal L},{\cal W}_{7}\right\} 
$. If ${\cal T}\in \left\{ {\cal C},{\cal W}_{0},{\cal W}_{8}\right\} $ none
of the condition is satisfied.
\end{rem-new}

\begin{th}
An Einstein manifold is globally $\varphi $-${\cal T}$-symmetric iff it is
globally $\varphi $-symmetric and $a_{0}\not=0$.
\end{th}

\noindent {\bf Proof.} By using (\ref{eq-GCT-1}) and (\ref{eq-phi-T-sym}),
we have the result.

\begin{th}
\label{thm-1} Let $M$ be a $3$-dimensional $\left( \varepsilon \right) $%
-para Sasakian manifold. Then $M$ is locally $\varphi $-${\cal T}$-symmetric
manifold if and only if the scalar curvature $r$ is constant.
\end{th}

\noindent {\bf Proof.} Let $M$ be a $3$-dimensional $\left( \varepsilon
\right) $-para Sasakian manifold. Differentiate covariantly on both sides of
(\ref{eq-S-5}), we have 
\begin{eqnarray}
(\nabla _{W}{\cal T})(X,Y)Z &=&\frac{\nabla _{W}r}{2}\left(
a_{0}+a_{1}+a_{4}+2a_{7}\right) g(Y,Z)X-\frac{\nabla _{W}r}{2}\left(
a_{0}-a_{2}-a_{5}+2a_{7}\right) g(X,Z)Y  \nonumber \\
&&+\frac{\nabla _{W}r}{2}\left( a_{3}+a_{6}\right) g(X,Y)Z-\frac{\nabla _{W}r%
}{2}a_{3}\eta (X)\eta (Y)Z-\left( \frac{\varepsilon r}{2}+3\right)
a_{3}\left( \nabla _{W}\eta \right) (X)\eta (Y)Z  \nonumber \\
&&-\left( \frac{\varepsilon r}{2}+3\right) a_{3}\eta (X)\left( \nabla
_{W}\eta \right) (Y)Z-\frac{\nabla _{W}r}{2}\left( a_{0}+a_{1}\right) \eta
(Y)\eta (Z)X  \nonumber \\
&&-\left( \frac{\varepsilon r}{2}+3\right) \left( a_{0}+a_{1}\right) \left(
\nabla _{W}\eta \right) (Y)\eta (Z)X-\left( \frac{\varepsilon r}{2}+3\right)
\left( a_{0}+a_{1}\right) \eta (Y)\left( \nabla _{W}\eta \right) (Z)X 
\nonumber \\
&&+\frac{\nabla _{W}r}{2}\left( a_{0}-a_{2}\right) \eta (X)\eta (Z)Y+\left( 
\frac{\varepsilon r}{2}+3\right) \left( a_{0}-a_{2}\right) \left( \nabla
_{W}\eta \right) (X)\eta (Z)Y  \nonumber \\
&&+\left( \frac{\varepsilon r}{2}+3\right) \left( a_{0}-a_{2}\right) \eta
(X)\left( \nabla _{W}\eta \right) (Z)Y+\frac{\nabla _{W}r}{2}\left(
a_{0}-a_{5}\right) g(X,Z)\eta (Y)\xi  \nonumber \\
&&+\left( \frac{r}{2}+3\varepsilon \right) \left( a_{0}-a_{5}\right)
g(X,Z)\left( \nabla _{W}\eta \right) (Y)\xi +\left( \frac{r}{2}+3\varepsilon
\right) \left( a_{0}-a_{5}\right) g(X,Z)\eta (Y)\nabla _{W}\xi  \nonumber \\
&&-\left( \frac{r}{2}+3\varepsilon \right) \left( a_{0}+a_{4}\right)
g(Y,Z)\left( \nabla _{W}\eta \right) (X)\xi -\left( \frac{r}{2}+3\varepsilon
\right) \left( a_{0}+a_{4}\right) g(Y,Z)\eta (X)\nabla _{W}\xi  \nonumber \\
&&-\left( \frac{r}{2}+3\varepsilon \right) a_{6}g(X,Y)\eta (Z)\nabla _{W}\xi
-\frac{\nabla _{W}r}{2}\left( a_{0}+a_{4}\right) g(Y,Z)\eta (X)\xi  \nonumber
\\
&&-\,\frac{\nabla _{W}r}{2}a_{6}g(X,Y)\eta (Z)\xi -\left( \frac{r}{2}%
+3\varepsilon \right) a_{6}g(X,Y)\left( \nabla _{W}\eta \right) (Z)\xi .
\label{eq-phi-T-sym-1}
\end{eqnarray}%
Applying $\varphi ^{2}$ on both sides of (\ref{eq-phi-T-sym-1}), we have 
\begin{eqnarray}
\varphi ^{2}(\nabla _{W}{\cal T})(X,Y)Z &=&\frac{\nabla _{W}r}{2}\left(
a_{0}+a_{1}+a_{4}+2a_{7}\right) g(Y,Z)(X-\eta (X)\xi )  \nonumber \\
&&-\frac{\nabla _{W}r}{2}\left( a_{0}-a_{2}-a_{5}+2a_{7}\right)
g(X,Z)(Y-\eta (Y)\xi )  \nonumber \\
&&+\frac{\nabla _{W}r}{2}\left( a_{3}+a_{6}\right) g(X,Y)(Z-\eta (Z)\xi )-%
\frac{\nabla _{W}r}{2}a_{3}\eta (X)\eta (Y)(Z-\eta (Z)\xi )  \nonumber \\
&&-\left( \frac{\varepsilon r}{2}+3\right) a_{3}\left( \nabla _{W}\eta
\right) (X)\eta (Y)(Z-\eta (Z)\xi )  \nonumber \\
&&-\left( \frac{\varepsilon r}{2}+3\right) a_{3}\eta (X)\left( \nabla
_{W}\eta \right) (Y)(Z-\eta (Z)\xi )  \nonumber \\
&&-\frac{\nabla _{W}r}{2}\left( a_{0}+a_{1}\right) \eta (Y)\eta (Z)(X-\eta
(X)\xi )  \nonumber \\
&&-\left( \frac{\varepsilon r}{2}+3\right) \left( a_{0}+a_{1}\right) \left(
\nabla _{W}\eta \right) (Y)\eta (Z)(X-\eta (X)\xi )  \nonumber \\
&&-\left( \frac{\varepsilon r}{2}+3\right) \left( a_{0}+a_{1}\right) \eta
(Y)\left( \nabla _{W}\eta \right) (Z)(X-\eta (X)\xi )  \nonumber \\
&&+\frac{\nabla _{W}r}{2}\left( a_{0}-a_{2}\right) \eta (X)\eta (Z)(Y-\eta
(Y)\xi )  \nonumber \\
&&+\left( \frac{\varepsilon r}{2}+3\right) \left( a_{0}-a_{2}\right) \left(
\nabla _{W}\eta \right) (X)\eta (Z)(Y-\eta (Y)\xi )  \nonumber \\
&&+\left( \frac{\varepsilon r}{2}+3\right) \left( a_{0}-a_{2}\right) \eta
(X)\left( \nabla _{W}\eta \right) (Z)(Y-\eta (Y)\xi )  \nonumber \\
&&+\left( \frac{r}{2}+3\varepsilon \right) \left( a_{0}-a_{5}\right)
g(X,Z)\eta (Y)\varphi ^{2}\nabla _{W}\xi  \nonumber \\
&&-\left( \frac{r}{2}+3\varepsilon \right) a_{6}g(X,Y)\eta (Z)\varphi
^{2}\nabla _{W}\xi  \nonumber \\
&&-\left( \frac{r}{2}+3\varepsilon \right) \left( a_{0}+a_{4}\right)
g(Y,Z)\eta (X)\varphi ^{2}\nabla _{W}\xi .  \label{eq-phi-T-sym-2+}
\end{eqnarray}%
Using the fact that $X$, $Y$, $Z$ are horizontal vector fields in (\ref%
{eq-phi-T-sym-2+}), we get 
\begin{eqnarray}
\varphi ^{2}(\nabla _{W}{\cal T})(X,Y)Z &=&\frac{\nabla _{W}r}{2}\left(
a_{0}+a_{1}+a_{4}+2a_{7}\right) g(Y,Z)X  \nonumber \\
&&-\frac{\nabla _{W}r}{2}\left( a_{0}-a_{2}-a_{5}+2a_{7}\right) g(X,Z)Y 
\nonumber \\
&&+\frac{\nabla _{W}r}{2}\left( a_{3}+a_{6}\right) g(X,Y)Z.
\label{eq-phi-T-sym-2.}
\end{eqnarray}%
If one of them $a_{0}+a_{1}+a_{4}+2a_{7}$, $a_{0}-a_{2}-a_{5}+2a_{7}$ and $%
a_{3}+a_{6}$ is not equal to zero, then by using (\ref{eq-phi-T-sym}), we
get the result. $\blacksquare $

\begin{rem-new}
One of them $a_{0}+a_{1}+a_{4}+2a_{7}$, $a_{0}-a_{2}-a_{5}+2a_{7}$ and $%
a_{3}+a_{6}$ is not equal to zero, for all the known curvature tensors.
\end{rem-new}

\section{$\protect\eta $-parallel Ricci tensor}

\begin{defn-new}
The Ricci tensor $S$ of an $\left( \varepsilon \right) $-para-Sasakian
manifold is called $\eta $-parallel if it satisfies 
\[
\left( \nabla _{X}S\right) \left( \varphi Y,\varphi Z\right) =0 
\]%
for all vector fields $X$, $Y$ and $Z$.
\end{defn-new}

\begin{th}
\label{thm-2} In a $3$-dimensional $\left( \varepsilon \right) $-para
Sasakian manifold with $\eta $-parallel Ricci tensor, the scalar curvature $%
r $ is constant.
\end{th}

\noindent Proof. By equation (\ref{eq-S-3}), we get 
\begin{equation}
S(\varphi Y,\varphi Z)=\left( \frac{r}{2}+\varepsilon \right) \left(
g(Y,Z)-\varepsilon \eta (Y)\eta (Z)\right)  \label{eq-S-311}
\end{equation}

Differentiating (\ref{eq-S-311}) covariantly with respect to $X$, we get 
\[
\left( \nabla _{X}S\right) (\varphi Y,\varphi Z)=\frac{\nabla _{X}r}{2}%
\left( g(Y,Z)-\varepsilon \eta (Y)\eta (Z)\right) -\varepsilon \left( \frac{r%
}{2}+\varepsilon \right) \left( \left( \nabla _{X}\eta \right) (Y)\eta
(Z)+\eta (Y)\left( \nabla _{X}\eta \right) (Z)\right) 
\]%
Suppose the Ricci tensor is $\eta $-parallel. Then from the above, we obtain 
\begin{equation}
\frac{\nabla _{X}r}{2}\left( g(Y,Z)-\varepsilon \eta (Y)\eta (Z)\right)
=\varepsilon \left( \frac{r}{2}+\varepsilon \right) \left( \left( \nabla
_{X}\eta \right) (Y)\eta (Z)+\eta (Y)\left( \nabla _{X}\eta \right)
(Z)\right)  \label{eq-S-3111}
\end{equation}%
Let $\{e_{i}\}$, $i=1,2,3$ be the orthonormal basis of tangent space at each
point of the manifold. Taking $Y=e_{i}=Z$ in (\ref{eq-S-3111}), we have $%
\nabla _{X}r=0$. Hence scalar curvature $r$ is constant. \medskip

\noindent From Theorems \ref{thm-1} and \ref{thm-2}, we can state the
following:

\begin{cor}
A $3$-dimensional $\left( \varepsilon \right) $-para Sasakian manifold with $%
\eta $-parallel Ricci tensor is locally $\varphi $-${\cal T}$-symmetric.
\end{cor}

\section{Example of a locally $\protect\varphi $-${\cal T}$-symmetric $%
\left( \protect\varepsilon \right) $-para Sasakian manifold of dimension $3$}

Consider the $3$-dimensional manifold $M=\left\{ (x,y,z)\in 
\mathbb{R}
^{3},z\neq 0\right\} $, where $(x,y,z)$ are the standard coordinates of $%
\mathbb{R}
^{3}$. The vector fields 
\[
e_{1}=z\frac{\partial }{\partial x},\qquad e_{2}=z\frac{\partial }{\partial y%
},\qquad e_{3}=-\;z\frac{\partial }{\partial z}
\]%
are linearly independent at each point of $M$. Let $g$ be the
semi-Riemannian metric defined by 
\[
\begin{array}{c}
g(e_{1},e_{3})=0,\qquad g(e_{1},e_{2})=0,\qquad g(e_{2},e_{3})=0, \\ 
g(e_{1},e_{1})=1,\qquad g(e_{2},e_{2})=1,\qquad g(e_{3},e_{3})=\varepsilon ,%
\end{array}%
\]%
where $\varepsilon =\pm 1$. Let $\eta $ be the $1$-form defined by $\eta
(Z)=\varepsilon g(Z,e_{3})$ for any $Z\in TM$. Let $\varphi $ be the $(1,1)$%
-tensor field defined by 
\[
\varphi e_{1}=\varepsilon e_{1},\qquad \varphi e_{2}=\varepsilon
e_{2},\qquad \varphi e_{3}=0.
\]%
Using the linearity of $\varphi $ and $g$, we have 
\[
\varphi ^{2}X=X-\eta (X)e_{3},
\]%
\[
\eta (e_{3})=1,
\]%
\[
g(\varphi X,\varphi Y)=g(X,Y)-\varepsilon \eta (X)\eta (Y),
\]%
\[
g(X,e_{3})=\varepsilon \eta (X),
\]%
\[
(\nabla _{X}\varphi )Y=-g(\varphi X,\varphi Y)e_{3}-\varepsilon \eta
(Y)\varphi ^{2}X,
\]%
for any $X,Y\in TM$. Then for $\xi =e_{3}$, the structure $\left( \varphi
,\xi ,\eta ,g,\varepsilon \right) $ defines an $\left( \varepsilon \right) $%
-para Sasakian structure on $M$. Let $\nabla $ be the Levi-Civita connection
with respect to the metric $g$. Then we have 
\[
\left[ e_{1},e_{2}\right] =0,\qquad \left[ e_{1},e_{3}\right] =e_{1},\qquad %
\left[ e_{1},e_{2}\right] =e_{2}.
\]%
The Koszul's formula for the Riemannian connection $\nabla $ of the metric $g
$ is given by 
\begin{eqnarray*}
2g(\nabla _{X}Y,Z) &=&Xg(Y,Z)+Yg(Z,X)-Zg(X,Y) \\
&&-g(X,[Y,Z])-g(Y,[X,Z])+g(Z,[X,Y]).
\end{eqnarray*}%
By using Koszul's formula, we have 
\[
\begin{array}{ccc}
\nabla _{e_{1}}e_{1}=-\varepsilon e_{3}, & \nabla _{e_{2}}e_{1}=0, & \nabla
_{e_{3}}e_{1}=-e_{1}, \\ 
\nabla _{e_{1}}e_{2}=0, & \nabla _{e_{2}}e_{2}=-\varepsilon e_{3}, & \nabla
_{e_{3}}e_{2}=-e_{2}, \\ 
\nabla _{e_{1}}e_{3}=e_{1}, & \nabla _{e_{2}}e_{3}=e_{2}, & \nabla
_{e_{3}}e_{3}=0.%
\end{array}%
\]%
From the above results, it is easy to check that equations (\ref{eq-1-1}), (%
\ref{eq-metric-1}), (\ref{eq-metric-2}), (\ref{eq-metric-3}), (\ref%
{eq-g(xi,xi)}) and (\ref{eq-eps-PS-def-1}) hold. Hence the manifold is an $%
(\varepsilon )$-para Sasakian manifold. \ 

Using the above results, it is easy to find out the following results 
\[
\begin{array}{ccc}
R(e_{1},e_{2})e_{1}=\varepsilon e_{2}, & R(e_{2},e_{3})e_{1}=0, & 
R(e_{1},e_{3})e_{1}=2\varepsilon e_{3}, \\ 
R(e_{1},e_{2})e_{2}=-\varepsilon e_{1}, & R(e_{2},e_{3})e_{2}=2\varepsilon
e_{3}, & R(e_{1},e_{3})e_{2}=0, \\ 
R(e_{1},e_{2})e_{3}=0, & R(e_{2},e_{3})e_{3}=0, & R(e_{1},e_{3})e_{3}=0.%
\end{array}%
\]%
Then 
\[
\begin{array}{ccc}
S(e_{1},e_{1})=-\left( \varepsilon +2\right) , & S(e_{2},e_{2})=-\left(
\varepsilon +2\right) , & S(e_{3},e_{3})=0,%
\end{array}%
\]%
and 
\[
r=-2\left( \varepsilon +2\right) . 
\]%
Hence the scalar curvature $r$ is constant. From Theorem \ref{thm-1}, $M$ is
a $3$-dimensional locally $\varphi $-${\cal T}$-symmetric $(\varepsilon )$
-para Sasakian manifold.

\medskip

\noindent Department of Mathematics\newline
Dr. Hari Singh Gour University\newline
Sagar-470003, Madhya Pradesh, India\newline
punam\_2101@yahoo.co.in

\end{document}